\documentclass[a4paper,11pt]{elsarticle}
\usepackage{multirow}
\usepackage[fleqn]{amsmath}
\usepackage{graphicx}
\usepackage{caption}
\usepackage{dsfont}
\usepackage[colorlinks,citecolor=blue,dvipdfm]{hyperref}
\usepackage{amsfonts}
\usepackage{subfigure}
\usepackage{lineno,hyperref}
\modulolinenumbers[5]
\usepackage{CJK}
\usepackage{indentfirst}
\usepackage{amsmath}
\usepackage{geometry}

\journal{}

\begin{document}

\begin{frontmatter}


\title{Stochastic Navier-Stokes equations with Caputo derivative driven by fractional noises}

\author[author1]{Guang-an Zou}
\author[author1]{Guangying Lv}
\author[author2]{Jiang-Lun Wu\corref{cor1}}
\cortext[cor1]{Corresponding author}
\ead{j.l.wu@swansea.ac.uk}

\address[author1]{School of Mathematics and Statistics, Henan University, Kaifeng 475004, P. R. China}
\address[author2]{Department of Mathematics, Swansea University, Swansea SA2 8PP, United Kingdom}

\begin{abstract}

In this paper, we consider the extended stochastic Navier-Stokes equations with Caputo derivative driven by fractional
Brownian motion. We firstly derive the pathwise spatial and temporal regularity of the generalized Ornstein-Uhlenbeck process. Then we discuss the existence, uniqueness, and H\"{o}lder regularity of  mild solutions to the given problem under certain sufficient conditions, which depend on the fractional order $\alpha$ and Hurst parameter $H$. The results obtained in this study improve some results in existing literature.
\end{abstract}


\end{frontmatter}


\section{Introduction}

Stochastic Navier-Stokes equations (SNSEs) are widely regarded as one of the most fascinating problems in fluid mechanics, in particular, stochasticity could even lead to a better understanding of physical phenomenon and mechanisms of turbulence in fluids. Furthermore, the presence of noises could give rise to some statistical features and important phenomena, for example, a unique invariant measure and ergodic behavior for the SNSEs driven by degenerate noise have been established, which can not be found in deterministic Navier-Stokes equations \cite{Flandoli-Maslowski-1995,Hairer-2006}. Since the seminal work of Bensoussan and Temam \cite{Bensoussan-1973}, the SNSEs have been intensively investigated in the literature. The existence and uniqueness of solutions for the SNSEs with multiplicative Gaussian noise were proved in \cite{Da-2002,Mikulevicius-2005,Taniguchi-2011}. The large deviation principle for SNSEs with multiplicative noise had been established in \cite{Wang-2015,Xu-2009}. The study of random attractors of SNSEs can be found in \cite{Brzezniak-2013, Flandoli-1995}, just mention a few.

On the other hand, fractional calculus has gained considerable popularity during the past decades owing to its demonstrated to describe physical systems possessing long-term memory and long-range spatial interactions, which play an important roles in diverse areas of science and engineering. Some theoretical analysis and experimental data have shown that fractional derivative can be recognized as one of the best tools to model anomalous diffusion processes \cite{Podlubny-1999,Srivastava-2006,Zhou-Wang-2016}. Consequently, the generalized Navier-Stokes equations with fractional derivative can be introduced to simulate anomalous diffusion in fractal media \cite{Momani-2006,Zhou-2017}. Recently, time-fractional Navier-Stokes equations have been initiated from the perspective of both analytical and numerical solutions, see \cite{De-2015,Ganji-2010,Kumar-2015,Li-2016,Zhou-Peng-2017,Zou-Zhou-2017} for more details.
We would like to emphasize that it is natural and also important to study the generalized SNSEs with time-fractional derivatives, which might be useful to model reasonably the phenomenon of anomalous diffusion with intrinsic random effects.

In this paper, we are concerned with the following generalized stochastic Navier-Stokes equations with time-fractional derivative on a finite time interval $[0,T]$ driven by fractional noise, defined on a domain $D\subset\mathbb{R}^d,d\ge1,$ with regular boundary $\partial D$
\begin{align*}
 ^{C}D_{t}^{\alpha}u=\nu\Delta u-(u\cdot\nabla)u-\nabla p+f(u)+\dot{B}^{H},~x\in D,~t>0, \tag{1.1}
\end{align*}
with the incompressibility condition:
\begin{align*}
\mathrm{div}  u=0,~x\in D,t\geq0, \tag{1.2}
\end{align*}
subject to the initial condition:
\begin{align*}
 u(x,0)=u_{0}(x),~x\in D,t=0, \tag{1.3}
\end{align*}
and the Dirichlet boundary conditions:
\begin{align*}
 u(x,t)=0,~x\in \partial D,t\geq0, \tag{1.4}
\end{align*}
in which $u=u(x,t)$ represents the velocity field of the fluid; $\nu>0$ is the viscosity coefficient; $p=p(x,t)$ is the associated pressure field; $f(u)$ stands for the deterministic external forces; The term $\dot{B}^{H}=\frac{d}{dt}B^{H}(t)$, and $B^{H}(t)$ represents a cylindrical fractional Brownian motion (fBm) with Hurst parameter $H\in(0,1)$ describes a state dependent random noise. Here, $^{C}D_{t}^{\alpha}$ denotes the Caputo-type derivative of order $\alpha$ ($0<\alpha<1$) for a function $u(x,t)$ with respect to time $t$ defined by
\begin{align*}
 ^{C}D_{t}^{\alpha}u(x,t)=\begin{cases}
\frac{1}{\Gamma(1-\alpha)}\int_{0}^{t}\frac{\partial u(x,s)}{\partial s}\frac{ds}{(t-s)^{\alpha}},~0<\alpha<1,\\
\frac{\partial u(x,t)}{\partial t}, ~~~~~~~~~~~~~~~~~~~~~~ \alpha =1,\\
\end{cases} \tag{1.5}
\end{align*}
where $\Gamma(\cdot)$ stands for the gamma function $\Gamma(x)=\int_{0}^{\infty}t^{x-1}e^{-t}dt$. Define the Stokes operator subject to the no-slip homogeneous Dirichlet boundary condition (1.4) as the formula
\begin{align*}
Au:=-\nu P_{H}\Delta u,
\end{align*}
where $P_{H}$ is the Helmhotz-Hodge projection operator, we also define the nonlinear operator $B$ as
\begin{align*}
B(u,v)=-P_{H}(u\cdot\nabla)v,
\end{align*}
with a slight abuse of notation $B(u):=B(u,u)$.

By applying the Helmholtz-Hodge operator $P_{H}$ to each term of time-fractional SNSEs, we can rewrite the Eqs.(1.1)-(1.4) as follows in the abstract form:
\begin{align*}
\begin{cases}
^{C}D_{t}^{\alpha}u=-Au+B(u)+f(u)+\dot{B}^{H},t>0,\\
u(0)=u_{0},
\end{cases} \tag{1.6}
\end{align*}
where we shall also use the same notations $f(u)$ instead of $P_{H}f$, and the solutions of problem (1.6) is also the solutions of Eqs.(1.1)-(1.4).

Note that the totality of fractional Brownian motions (fBms) form a subclass of Gaussian processes, which are positively correlated for $H\in(1/2,1)$ and negatively correlated for $H\in(0,1/2)$ while $H=1/2$ are standard Gaussian processes. So it is interesting to consider the stochastic differential equations with fBm, and the subject of stochastic calculus with respect to fBm has attracted much attentions \cite{Biagini-2008,Duncan-2009,Jiang-2012,Mishura-2008}. In recent years, the existence and uniqueness of solutions for stochastic Burgers equations with fBm have been examined in \cite{Jiang-2012,Wang-Zeng-2010}. In addition, Zou and Wang investigated the time-fractional stochastic Burgers equation with standard Brownian motion \cite{Zou-2017}. However, to the best of our knowledge, the study of time-fractional SNSEs with fBm has not been addressed yet, which is indeed a fascinating and more interesting (and also practical) problem. The objective of the present paper is to establish the existence and uniqueness of mild solutions by Banach fixed point theorem and Mainardi's Wright-type function, the  key and difficulty of problems are how to deal with stochastic convolution. We also prove the H\"{o}lder regularity of mild solutions to time-fractional SNSEs. Our consideration extends and improves the existing results carried out in previous studies \cite{De-2015,Flandoli-Maslowski-1995,Mikulevicius-2005,Wang-Zeng-2010,Zhou-2017}.

The rest of the paper is organized as follows. In the next section, we introduce several notions and we give certain preliminaries needed in our later analysis. In Section 3, we establish the pathwise spatial and temporal regularity of the generalized Ornstein-Uhlenbeck process. In Section 4, we show the existence and uniqueness of mild solutions to time-fractional SNSEs. We end our paper by proving the H\"{o}lder regularity of the mild solution.

\section{Notations and preliminaries}

In this section, we give some notions and certain important preliminaries, which will be used in the subsequent discussions.

Let $(\Omega,\mathcal{F},\mathds{P},\{\mathcal{F}_{t}\}_{t\geq0})$ be a filtered probability space with a normal filtration $\{\mathcal{F}_{t}\}_{t\geq0}$. We assume that the operator $A$ is self-adjoint and there exist the eigenvectors $e_{k}$ corresponding to eigenvalues $\gamma_{k}$ such that
\begin{align*}
Ae_{k}=\gamma_{k}e_{k},e_{k}=\sqrt{2}\sin(k\pi),\gamma_{k}=\pi^{2}k^{2},k\in N^{+}.
\end{align*}

For any $\sigma>0$, $A^{\frac{\sigma}{2}}e_{k}=\gamma_{k}^{\frac{\sigma}{2}}e_{k}, k=1,2,\ldots$, and let $\dot{H}^{\sigma}$ be the domain of the fractional power defined by
\begin{align*}
 \dot{H}^{\sigma}=\mathcal{D}(A^{\frac{\sigma}{2}})=\{v\in L^{2}(D),s.t.~\|v\|_{ \dot{H}^{\sigma}}^{2}=\sum\limits_{k=1}^{\infty}\gamma_{k}^{\frac{\sigma}{2}}v_{k}^{2}<\infty\},
\end{align*}
where $v_{k}=(v,e_{k})$ and the norm $\|v\|_{ \dot{H}^{\sigma}}=\|A^{\frac{\sigma}{2}}v\|$. Let $L^{2}(\Omega,H)$ be a Hilbert space of $H$-valued random variables equipped with the inner product $\mathbb{E}(\cdot,\cdot)$ and norm $\mathbb{E}\|\cdot\|$, it is given by
\begin{align*}
L_{2}(\Omega,H)=\{\chi:\mathbb{E}\|\chi\|_{H}^{2}=\int_{\Omega}\|\chi(\omega)\|_{H}^{2}d\mathds{P}(\omega)<\infty,\omega\in \Omega\}.
\end{align*}

\textbf{Definition 2.1.} For $H\in(0,1)$, a continuous centered Gaussian process $\{\beta^{H}(t),t\in[0,\infty)\}$ with covariance function
\begin{align*}
R_{H}(t,s)=\mathbb{E}[\beta^{H}(t)\beta^{H}(s)]=\frac{1}{2}(t^{2H}+s^{2H}-|t-s|^{2H}),~t,s\in [0,\infty)
\end{align*}
is called a one-dimensional fractional Brownian motion (fBm), and $H$ is the Hurst parameter. In particular when $H=\frac{1}{2}$, $\beta^{H}(t)$ represents a standard Brownian motion.

Now let us introduce the Wiener integral with respect to the fBm. To begin with, we represent  $\beta^{H}(t)$ as following (see \cite{Biagini-2008})
\begin{align*}
\beta^{H}(t)=\int_{0}^{t}K_{H}(t,s)dW(s),
\end{align*}
where $W=\{W(t),t\in[0,T]\}$ is a Wiener process on the space $(\Omega,\mathcal{F},\mathds{P},\{\mathcal{F}_{t}\}_{t\geq0})$ and the kernel $K_{H}(t,s), 0\le s< t\le T$, is given by
\begin{align*}
K_{H}(t,s):=c_{H}(t-s)^{H-\frac{1}{2}}+c_{H}(\frac{1}{2}-H)\int_{s}^{t}(u-s)^{H-\frac{3}{2}}(1-(\frac{s}{u})^{\frac{1}{2}-H})du, \tag{2.1}
\end{align*}
for $0<H<\frac{1}{2}$ and $c_{H}=(\frac{2H\Gamma(\frac{3}{2}-H)}{\Gamma(H+\frac{1}{2})\Gamma(2-2H)})^{\frac{1}{2}}$ is a constant. When $\frac{1}{2}<H<1$, there holds
\begin{align*}
K_{H}(t,s)=c_{H}(H-\frac{1}{2})s^{\frac{1}{2}-H}\int_{s}^{t}(u-s)^{H-\frac{3}{2}}u^{H-\frac{1}{2}}du.\tag{2.2}
\end{align*}

It is easy to verify that
\begin{align*}
\frac{\partial K_{H}}{\partial t}(t,s)=c_{H}(H-\frac{1}{2})(\frac{s}{t})^{\frac{1}{2}-H}(t-s)^{H-\frac{3}{2}}. \tag{2.3}
\end{align*}

We denote by $\mathcal{H}$ the reproducing kernel Hilbert space of the fBm. Let $K_{\tau}^{*}:\mathcal{H}\rightarrow L^{2}([0,T])$ be the linear map given by
\begin{align*}
(K_{\tau}^{*}\psi)(s)=\varphi(s)K_{H}(\tau,s)+\int_{s}^{\tau}(\psi(t)-\psi(s))\frac{\partial K_{H}}{\partial t}(t,s)dt  \tag{2.4}
\end{align*}
for $0<H<\frac{1}{2}$, and if $\frac{1}{2}<H<1$, we denote
\begin{align*}
(K_{\tau}^{*}\psi)(s)=\int_{s}^{\tau}\psi(t)\frac{\partial K_{H}}{\partial t}(t,s)dt. \tag{2.5}
\end{align*}

We refer the reader to \cite{Mishura-2008} for the proof of the fact that $K_{\tau}^{*}$ is an isometry between $\mathcal{H}$ and $L^{2}([0,T])$. Moreover, for any $\psi\in \mathcal{H}$, we have the following relation between the Wiener integral with respect to fBm and the It\^{o} integral with respect to Wiener process
\begin{align*}
\int_{0}^{t}\psi(s)d\beta^{H}(s)=\int_{0}^{t}(K_{\tau}^{*}\psi)(s)dW(s),~t\in [0,T].
\end{align*}

Generally, following the standard approach for $H=\frac{1}{2}$, we consider $Q$-Wiener process with linear bounded covariance operator $Q$ such that $\mathrm{Tr} (Q)<\infty$. Furthermore, there exists the eigenvalues $\lambda_{n}$ and corresponding eigenfunctions $e_{k}$ satisfying $Q e_{k}=\lambda_{n}e_{k},k=1,2,\ldots$, then we define the infinite dimensional fBm with covariance $Q$ as
\begin{align*}
B^{H}(t):=\sum\limits_{k=1}^{\infty}\lambda^{1/2}_{k}e_{k}\beta_{k}^{H}(t),
\end{align*}
where $\beta_{k}^{H}$ are real-valued independent fBm's. In order to define Wiener integrals with repect to $Q$-fBm, we introduce $\mathcal{L}_{2}^{0}:=\mathcal{L}_{2}^{0}(Y,X)$ of all $Q$-Hilbert-Schmidt operators $\psi:Y\rightarrow X$, where $Y$ and $X$ are two real separable Hilbert spaces. We associate the $Q$-Hilbert-Schmidt operators $\psi$ with the norm
\begin{align*}
\|\varphi\|_{\mathcal{L}_{2}^{0}}^{2}=\sum\limits_{k=1}^{\infty}\|\lambda^{1/2}_{k}\psi e_{k}\|^{2}<\infty.
\end{align*}

As a consequence, for $\psi\in \mathcal{L}_{2}^{0}(Y,X)$, the Wiener integral of $\psi$ with respect to $B^{H}(t)$ is defined by
\begin{align*}
\int_{0}^{t}\psi(s)dB^{H}(s)=\sum\limits_{k=1}^{\infty}\int_{0}^{t}\lambda^{1/2}_{k}\psi(s)e_{k}d\beta_{k}^{H}(s)=\sum\limits_{k=1}^{\infty}\int_{0}^{t}\lambda^{1/2}_{k}(K_{\tau}^{*}\psi e_{k})(s)d\beta_{k}(s), \tag{2.6}
\end{align*}
where $\beta_{k}$ is the standard Brownian motion.

\textbf{Definition 2.2.}  An $\mathcal{F}_{t}$-adapted stochastic process $(u(t),t\in[0,T])$ is called a mild solution to (1.6) if
the following integral equation is satisfied
\begin{align*}
 u(t)&=E_{\alpha}(t)u_{0}+\int_{0}^{t}(t-s)^{\alpha-1}E_{\alpha,\alpha}(t-s)[B(u(s))+f(u(s))]ds\\
 &~~~+\int_{0}^{t}(t-s)^{\alpha-1}E_{\alpha,\alpha}(t-s)dB^{H}(s), \tag{2.7}
\end{align*}
where the generalized Mittag-Leffler operators $E_{\alpha}(t)$ and $E_{\alpha,\alpha}(t)$ are defined, respectively, by
\begin{align*}
E_{\alpha}(t):=\int_{0}^{\infty}\xi_{\alpha}(\theta)T(t^{\alpha}\theta)d\theta,
\end{align*}
and
\begin{align*}
E_{\alpha,\alpha}(t):=\int_{0}^{\infty}\alpha\theta\xi_{\alpha}(\theta)T(t^{\alpha}\theta)d\theta,
\end{align*}
where $T(t)=e^{-tA},t\geq0$ is an analytic semigroup generated by the operator $-A$, and the Mainardi's Wright-type function with $\alpha\in (0,1)$ is given by
\begin{align*}
\xi_{\alpha}(\theta)=\sum_{k=0}^{\infty}\frac{(-1)^{k}\theta^{k}}{k!\Gamma(1-\alpha(1+k))}.
\end{align*}

Furthermore, for any $\alpha\in (0,1)$ and $-1<\nu<\infty$, it is not difficult to verity that
\begin{align*}
\xi_{\alpha}(\theta)\geq0 ~and~ \int_{0}^{\infty}\theta^{\nu}\xi_{\alpha}(\theta)d\theta=\frac{\Gamma(1+\nu)}{\Gamma(1+\alpha\nu)}, \tag{2.8}
\end{align*}
for all $\theta\geq0$. The derivation of mild solution (2.7) can refer to \cite{Zou-2017}.

The operators $\{E_{\alpha}(t)\}_{t\geq0}$ and $\{E_{\alpha,\alpha}(t)\}_{t\geq0}$ in (2.7) have the following properties \cite{Zou-2017}:

\textbf{Lemma 2.1.} For any $t>0$, $E_{\alpha}(t)$ and $E_{\alpha,\alpha}(t)$ are linear and bounded operators. Moreover, for $0<\alpha<1$ and $0\leq\nu<2$, there exists a constant $C>0$ such that
\begin{align*}
\|E_{\alpha}(t)\chi\|_{\dot{H}^{\nu}}\leq Ct^{-\frac{\alpha\nu}{2}}\|\chi\|,~\|E_{\alpha,\alpha}(t)\chi\|_{\dot{H}^{\nu}}\leq Ct^{-\frac{\alpha\nu}{2}}\|\chi\|.
\end{align*}

\textbf{Lemma 2.2.} For any $t>0$, the operators $E_{\alpha}(t)$ and $E_{\alpha,\alpha}(t)$ are strongly continuous. Moreover, for $0<\alpha<1$ and $0\leq\nu<2$ and $0\leq t_{1}< t_{2}\leq T$, there exists a constant $C>0$ such that
\begin{align*}
\|(E_{\alpha}(t_{2})-E_{\alpha}(t_{1}))\chi\|_{\dot{H}^{\nu}}\leq C(t_{2}-t_{1})^{\frac{\alpha\nu}{2}}\|\chi\|,
\end{align*}
and
\begin{align*}
\|(E_{\alpha,\alpha}(t_{2})-E_{\alpha,\alpha}(t_{1}))\chi\|_{\dot{H}^{\nu}}\leq C(t_{2}-t_{1})^{\frac{\alpha\nu}{2}}\|\chi\|.
\end{align*}

Throughout the paper, we assume that the mapping $f: \Omega\times H \rightarrow H$ satisfies
the following global Lipschitz and growth conditions
\begin{align*}
\|f(u)-f(v)\|^{2}\leq C\|u-v\|^{2},~\|f(u)\|^{2}\leq C(1+\|u\|^{2})\tag{2.9}
\end{align*}
for any $u,v\in H$.

\section{Regularity of the stochastic convolution}

In this section, we state and prove the basic properties of stochastic convolution. Firstly, we introduce the following generalized Ornstein-Uhlenbeck process
\begin{equation*}
Z(t):=\int_{0}^{t}(t-s)^{\alpha-1}E_{\alpha,\alpha}(t-s)dB^{H}(s). \tag{3.1}
\end{equation*}

Obviously, it is very important to establish the basic properties of the stochastic integrals (3.1)
in the study of the problem (1.6). For the sake of convenience, we introduce the following operator and show some properties.

\textbf{Lemma 3.1.} Let $\mathcal{S}_{\alpha}(t)=t^{\alpha-1}E_{\alpha,\alpha}(t)$, for $0\leq\nu< 2$ and $0<\alpha<1$, there exists a constant $C>0$ such that
\begin{align*}
\|\mathcal{S}_{\alpha}(t)\chi\|_{\dot{H}^{\nu}}\leq Ct^{\frac{(2-\nu)\alpha-2}{2}}\|\chi\|,~\|[\mathcal{S}_{\alpha}(t_{2})-\mathcal{S}_{\alpha}(t_{1})]\chi\|_{\dot{H}^{\nu}}\leq C(t_{2}-t_{1})^{\frac{ 2-(2-\nu)\alpha}{2}}\|\chi\|
\end{align*}
for any $0\leq t_{1}<t_{2}\leq T$.

\textbf{Proof.} By Lemma 2.1, we get
\begin{align*}
\|\mathcal{S}_{\alpha}(t)\chi\|_{\dot{H}^{\nu}}=\|t^{\alpha-1}E_{\alpha,\alpha}(t)\chi\|_{\dot{H}^{\nu}}\leq Ct^{\frac{(2-\nu)\alpha-2}{2}}\|\chi\|.
\end{align*}

Next, utilizing the property of semigroup $\|A^{\sigma}e^{-tA}\|\leq Ct^{-\sigma}$ for $\sigma\geq0$, we have
\begin{align*}
\|\frac{d}{dt}\mathcal{S}_{\alpha}(t)\chi \|_{\dot{H}^{\nu}}&=\|(\alpha-1)t^{\alpha-2}E_{\alpha,\alpha}(t)-\int_{0}^{\infty}\alpha^{2}t^{2\alpha-2}\theta^{2}\xi_{\alpha}(\theta)AT(t^{\alpha}\theta)d\theta\|_{\dot{H}^{\nu}}\\
&\leq (1-\alpha)t^{\alpha-2}\|E_{\alpha,\alpha}(t)\chi\|_{\dot{H}^{\nu}}+\int_{0}^{\infty}\alpha^{2}t^{2\alpha-2}\theta^{2}\xi_{\alpha}(\theta)\|A^{1+\frac{\nu}{2}}e^{-t^{\alpha}\theta A}\chi\|d\theta\\
&\leq C(1-\alpha)t^{\frac{(2-\nu)\alpha-4}{2}}\|\chi\|+\frac{\alpha^{2}\Gamma(2-\frac{\nu}{2})}{\Gamma(1+\alpha(1-\frac{\nu}{2}))}t^{\frac{(2-\nu)\alpha-4}{2}}\|\chi\|\\
&\leq Ct^{\frac{(2-\nu)\alpha-4}{2}}\|\chi\|.
\end{align*}
Hence, we have the following
\begin{align*}
\|[\mathcal{S}_{\alpha}(t_{2})-\mathcal{S}_{\alpha}(t_{1})]\chi\|_{\dot{H}^{\nu}}&=\|\int_{t_{1}}^{t_{2}}\frac{d}{dt}\mathcal{S}_{\alpha}(t)\chi dt\|_{\dot{H}^{\nu}}\\
&\leq\int_{t_{1}}^{t_{2}}Ct^{\frac{(2-\nu)\alpha-4}{2}}\|\chi\|dt\\
&=\frac{2C}{2-(2-\nu)\alpha}[t_{1}^{\frac{(2-\nu)\alpha-2}{2}}-t_{2}^{\frac{(2-\nu)\alpha-2}{2}}]\\
&\leq \frac{2C}{[2-(2-\nu)\alpha]T_{0}^{2-(2-\nu)\alpha}}(t_{2}-t_{1})^{\frac{2-(2-\nu)\alpha}{2}},
\end{align*}
where $0<T_{0}\leq t_{1}<t_{2}\leq T$, and we have used $t_{2}^{\omega}-t_{1}^{\omega}\leq C(t_{2}-t_{1})^{\omega}$ for $0\leq \omega\leq1$, in the above derivation.

In what follows, let us establish the pathwise spatial-temporal regularity of the stochastic convolution (3.1).

\textbf{Theorem 3.1.} For $0\leq\nu<2$ and $0<\alpha<1$, the generalized Ornstein-Uhlenbeck process $(Z(t))_{t\geq0}$ with the Hurst parameter $\frac{1}{4}<H<1$ is well defined. Moreover, there holds
\begin{align*}
\sup\limits_{t\in[0,T]}\mathbb{E}\|Z(t)\|_{\dot{H}^{\nu}}^{2}\leq C(H,Q)t^{\sigma}<\infty,~
0\leq t\leq T,
\end{align*}
where the index should satisfy $\sigma=\min\{(2-\nu)\alpha+4H-3,(2-\nu)\alpha+2H-1\}>0$.

\textbf{Proof.} Using the Wiener integral with respect to fBm and noticing the expression of $K_{t}^{*}$ and the properties of It\^{o} integral, for $0<H<\frac{1}{2}$, we get
\begin{align*}
\mathbb{E}\|Z(t)\|_{\dot{H}^{\nu}}^{2}&=\mathbb{E}\|\int_{0}^{t}(t-s)^{\alpha-1}E_{\alpha,\alpha}(t-s)dB^{H}(s)\|_{\dot{H}^{\nu}}^{2}\\
&=\sum\limits_{k=1}^{\infty}\mathbb{E}\|\int_{0}^{t}\lambda^{1/2}_{k}(K_{t}^{*}\mathcal{S}_{\alpha}(t-s)e_{k})(s)d\beta_{k}(s)\|_{\dot{H}^{\nu}}^{2}\\
&=\sum\limits_{k=1}^{\infty}\int_{0}^{t}\mathbb{E}\|\lambda^{1/2}_{k}(K_{t}^{*}\mathcal{S}_{\alpha}(t-s)e_{k})(s)\|_{\dot{H}^{\nu}}^{2}ds\\
&=\sum\limits_{k=1}^{\infty}\int_{0}^{t}\mathbb{E}\|\lambda^{1/2}_{k}\mathcal{S}_{\alpha}(t-s)K_{H}(t,s)e_{k}\\
&\hspace{2mm}+\int_{s}^{t}\lambda^{1/2}_{k}[\mathcal{S}_{\alpha}(t-r)-\mathcal{S}_{\alpha}(t-s)]\frac{\partial K_{H}}{\partial r}(r,s)e_{k}dr\|_{\dot{H}^{\nu}}^{2}ds\\
&\leq 2\sum\limits_{k=1}^{\infty}\int_{0}^{t}\mathbb{E}\|\lambda^{1/2}_{k}\mathcal{S}_{\alpha}(t-s)K_{H}(t,s)e_{k}\|_{\dot{H}^{\nu}}^{2}ds\\
&\hspace{2mm}+2\sum\limits_{k=1}^{\infty}\int_{0}^{t}\mathbb{E}\|\int_{s}^{t}\lambda^{1/2}_{k}[\mathcal{S}_{\alpha}(t-r)-\mathcal{S}_{\alpha}(t-s)]\frac{\partial K_{H}}{\partial r}(r,s)e_{k}dr\|_{\dot{H}^{\nu}}^{2}ds\\
&=:I_{1}+I_{2}. \tag{3.2}
\end{align*}

With the help of the following inequality (see \cite{Wang-Zeng-2010})
\begin{align*}
K_{H}(t,s)\leq C(H)(t-s)^{H-\frac{1}{2}}s^{H-\frac{1}{2}},
\end{align*}
and further combining Lemma 3.1 and the H\"{o}lder inequality, we obtain
\begin{align*}
I_{1}&=2\sum\limits_{k=1}^{\infty}\int_{0}^{t}\mathbb{E}\|\lambda^{1/2}_{k}\mathcal{S}_{\alpha}(t-s)K_{H}(t,s)e_{k}\|_{\dot{H}^{\nu}}^{2}ds\\
&\leq 2C(H)(\int_{0}^{t}(t-s)^{(2-\nu)\alpha+2H-3}s^{2H-1}\sum\limits_{k=1}^{\infty}\mathbb{E}\|\lambda^{1/2}_{k}e_{k}\|^{2}ds)\\
&\leq 2C(H)Tr(Q)(\int_{0}^{t}(t-s)^{2[(2-\nu)\alpha+2H-3]}ds)^{\frac{1}{2}}(\int_{0}^{t}s^{2(2H-1)}ds)^{\frac{1}{2}}\\
&\leq C(H,Q)t^{(2-\nu)\alpha+4H-3}, \tag{3.3}
\end{align*}
and on the other hand, utilizing the expression (2.3), we get
\begin{align*}
I_{2}&=2\sum\limits_{k=1}^{\infty}\int_{0}^{t}\mathbb{E}\|\int_{s}^{t}\lambda^{1/2}_{k}[\mathcal{S}_{\alpha}(t-r)-\mathcal{S}_{\alpha}(t-s)]\frac{\partial K_{H}}{\partial r}(r,s)e_{k}dr\|_{\dot{H}^{\nu}}^{2}ds\\
&\leq 2\sum\limits_{k=1}^{\infty}\int_{0}^{t}\mathbb{E}(\int_{s}^{t}\|[\mathcal{S}_{\alpha}(t-r)-\mathcal{S}_{\alpha}(t-s)]\frac{\partial K_{H}}{\partial r}(r,s)\|_{\dot{H}^{\nu}}^{2}dr)(\int_{s}^{t}\mathbb{E}\|\lambda^{1/2}_{k}e_{k}\|^{2}dr)ds\\
&\leq2C_{H}^{2}(H-\frac{1}{2})^{2}Tr(Q)\int_{0}^{t}(t-s)(\int_{s}^{t}|(s-r)^{\frac{(2-\nu)\alpha}{2}}(\frac{s}{r})^{\frac{1}{2}-H}(r-s)^{H-\frac{3}{2}}|^{2}dr)ds\\
&\leq C(H,Q)(\int_{0}^{t}(t-s)^{(2-\nu)\alpha+4H-3}s^{1-2H}ds)\\
&\leq C(H,Q)t^{(2-\nu)\alpha+2H-1}. \tag{3.4}
\end{align*}

When $\frac{1}{4}<H<\frac{1}{2}$ and $\sigma=\min\{(2-\nu)\alpha+4H-3,(2-\nu)\alpha+2H-1\}>0$, by combining (3.2)-(3.4), one can easily get that
\begin{align*}
\mathbb{E}\|Z(t)\|_{\dot{H}^{\nu}}^{2}\leq C(H,Q)t^{\nu}\leq C(H,Q)T^{\sigma}<\infty.
\end{align*}

Similarly, for $\frac{1}{2}<H<1$, one can derive that
\begin{align*}
&\mathbb{E}\|Z(t)\|_{\dot{H}^{\nu}}^{2}\\
&=\mathbb{E}\|\int_{0}^{t}(t-s)^{\alpha-1}E_{\alpha,\alpha}(t-s)dB^{H}(s)\|_{\dot{H}^{\nu}}^{2}\\
&=\sum\limits_{k=1}^{\infty}\int_{0}^{t}\mathbb{E}\|\lambda^{1/2}_{k}(K_{t}^{*}\mathcal{S}_{\alpha}(t-s)e_{k})(s)\|_{\dot{H}^{\nu}}^{2}ds\\
&=\sum\limits_{k=1}^{\infty}\int_{0}^{t}\mathbb{E}\|\int_{s}^{t}\lambda^{1/2}_{k}\mathcal{S}_{\alpha}(t-r)\frac{\partial K_{H}}{\partial r}(r,s)e_{k}dr\|_{\dot{H}^{\nu}}^{2}ds\\
&\leq C_{H}^{2}(H-\frac{1}{2})^{2}\int_{0}^{t}\mathbb{E}(\int_{s}^{t}\|\mathcal{S}_{\alpha}(t-r)(\frac{s}{r})^{\frac{1}{2}-H}(r-s)^{H-\frac{3}{2}}\|_{\dot{H}^{\nu}}^{2}dr)(\int_{s}^{t}\mathbb{E}\|\lambda^{1/2}_{k}e_{k}\|^{2}dr)ds\\
&\leq C(H,Q)(\int_{0}^{t}(t-s)^{(2-\nu)\alpha+4H-3}s^{1-2H}ds)\\
&\leq C(H,Q)t^{(2-\nu)\alpha+2H-1}.
\end{align*}

Thus, if $\frac{1}{2}<H<1$ and $(2-\nu)\alpha+2H-1>0$, one can directly obtain $\mathbb{E}\|Z(t)\|_{\dot{H}^{\nu}}^{2}<C(H,Q)T^{(2-\nu)\alpha+2H-1}<\infty$. When $H=\frac{1}{2}$, $B^{H}(t)$ is standard Brownian motion and it is easy to obtain $Z(t)$ is well defined. This completes the proof. $\square$

\textbf{Theorem 3.2.} For $0\leq\nu<2$ and $0<\alpha<1$, the stochastic process $(Z_{t})_{t\geq0}$ with $\frac{1}{4}<H<1$ is continuous and it satisfies
\begin{align*}
\mathbb{E}\|Z(t_{2})-Z(t_{1})\|_{\dot{H}^{\nu}}^{2}\leq C(H,Q)(t_{2}-t_{1})^{\gamma},~
0\leq t_{1}<t_{2}\leq T,
\end{align*}
where the index $\gamma=\min\{2-(2-\nu)\alpha,(2-\nu)\alpha+4H-3,(2-\nu)\alpha+2H-1\}>0$.

\textbf{Proof.} From (2.7), according to the relation between the Wiener integral and fBm, we have
\begin{align*}
Z(t_{2})-Z(t_{1})&=\int_{0}^{t_{2}}(t_{2}-s)^{\alpha-1}E_{\alpha,\alpha}(t_{2}-s)dB^{H}(s)-\int_{0}^{t_{1}}(t_{1}-s)^{\alpha-1}E_{\alpha,\alpha}(t_{1}-s)dB^{H}(s)\\
&=\int_{0}^{t_{1}}(\mathcal{S}_{\alpha}(t_{2}-s)-\mathcal{S}_{\alpha}(t_{1}-s))dB^{H}(s)+\int_{t_{1}}^{t_{2}}\mathcal{S}_{\alpha}(t_{2}-s)dB^{H}(s)\\
&=\sum\limits_{k=1}^{\infty}\int_{0}^{t_{1}}\lambda^{1/2}_{k}(K_{t}^{*}(\mathcal{S}_{\alpha}(t_{2}-s)-\mathcal{S}_{\alpha}(t_{1}-s))e_{k})(s)d\beta_{k}(s)\\
&\hspace{2mm}+\sum\limits_{k=1}^{\infty}\int_{t_{1}}^{t_{2}}\lambda^{1/2}_{k}(K_{t}^{*}\mathcal{S}_{\alpha}(t_{2}-s)e_{k})(s)d\beta_{k}(s)\\
&=:J_{1}+J_{2}. \tag{3.5}
\end{align*}

For the term $J_{1}$, we get
\begin{align*}
&\mathbb{E}\|J_{1}\|_{\dot{H}^{\sigma}}^{2}\\
&=\mathbb{E}\|\sum\limits_{k=1}^{\infty}\int_{0}^{t_{1}}\lambda^{1/2}_{k}(K_{t}^{*}(\mathcal{S}_{\alpha}(t_{2}-s)-\mathcal{S}_{\alpha}(t_{1}-s))e_{k})(s)d\beta_{k}(s)\|_{\dot{H}^{\nu}}^{2}\\
&=\sum\limits_{k=1}^{\infty}\int_{0}^{t_{1}}\mathbb{E}\|\lambda^{1/2}_{k}(K_{t}^{*}(\mathcal{S}_{\alpha}(t_{2}-s)-\mathcal{S}_{\alpha}(t_{1}-s))e_{k})(s)\|_{\dot{H}^{\nu}}^{2}ds\\
&=\sum\limits_{k=1}^{\infty}\int_{0}^{t_{1}}\mathbb{E}\|\lambda^{1/2}_{k}(\mathcal{S}_{\alpha}(t_{2}-s)-\mathcal{S}_{\alpha}(t_{1}-s))K_{H}(t,s)e_{k}\\
&\hspace{2mm}+\int_{s}^{t}\lambda^{1/2}_{k}[(\mathcal{S}_{\alpha}(t_{2}-r)-\mathcal{S}_{\alpha}(t_{1}-r))-(\mathcal{S}_{\alpha}(t_{2}-s)-\mathcal{S}_{\alpha}(t_{1}-s))]\frac{\partial K_{H}}{\partial r}(r,s)e_{k}dr\|_{\dot{H}^{\nu}}^{2}ds\\
&\leq 2 (t_{2}-t_{1})^{2-(2-\nu)\alpha}\int_{0}^{t}(\|K_{H}(t,s)\|^{2}\mathbb{E}\|\lambda^{1/2}_{k}e_{k}\|^{2}+2\mathbb{E}(\int_{s}^{t}\|\frac{\partial K_{H}}{\partial r}(r,s)\|^{2}dr)(\int_{s}^{t}\mathbb{E}\|\lambda^{1/2}_{k}e_{k}\|^{2}dr))ds\\
&\leq C(H)Tr(Q)(t_{2}-t_{1})^{2-(2-\nu)\alpha}\int_{0}^{t}[(t-s)^{2H-1}s^{2H-1}+(t-s)(\int_{s}^{t}(\frac{s}{r})^{1-2H}(r-s)^{2H-3}dr)]ds\\
&\leq C(H,Q)t^{2H}(t_{2}-t_{1})^{2-(2-\nu)\alpha}. \tag{3.6}
\end{align*}

Applying Lemma 3.1 and the H\"{o}lder inequality, we obtain
\begin{align*}
\mathbb{E}\|J_{2}\|_{\dot{H}^{\nu}}^{2}&=\mathbb{E}\|\sum\limits_{k=1}^{\infty}\int_{t_{1}}^{t_{2}}\lambda^{1/2}_{k}(K_{t}^{*}\mathcal{S}_{\alpha}(t_{2}-s)e_{k})(s)d\beta_{k}(s)\|_{\dot{H}^{\nu}}^{2}\\
&=\sum\limits_{k=1}^{\infty}\int_{t_{1}}^{t_{2}}\mathbb{E}\|\lambda^{1/2}_{k}(K_{t}^{*}\mathcal{S}_{\alpha}(t_{2}-s)e_{k})(s)\|_{\dot{H}^{\nu}}^{2}ds\\
&=\sum\limits_{k=1}^{\infty}\int_{t_{1}}^{t_{2}}\mathbb{E}\|\lambda^{1/2}_{k}\mathcal{S}_{\alpha}(t_{2}-s)K_{H}(t,s)e_{k}\\
&\hspace{2mm}+\int_{s}^{t}\lambda^{1/2}_{k}[\mathcal{S}_{\alpha}(t_{2}-r)-\mathcal{S}_{\alpha}(t_{2}-s)]\frac{\partial K_{H}}{\partial r}(r,s)e_{k}dr\|_{\dot{H}^{\nu}}^{2}ds\\
&\leq 2\int_{t_{1}}^{t_{2}}\|\mathcal{S}_{\alpha}(t_{2}-s)K_{H}(t,s)\|_{\dot{H}^{\nu}}^{2}\mathbb{E}\|\lambda^{1/2}_{k}e_{k}\|^{2}ds\\
&\hspace{2mm}+2\sum\limits_{k=1}^{\infty}\int_{t_{1}}^{t_{2}}\mathbb{E}\|\int_{s}^{t}\lambda^{1/2}_{k}[\mathcal{S}_{\alpha}(t_{2}-r)-\mathcal{S}_{\alpha}(t_{2}-s)]\frac{\partial K_{H}}{\partial r}(r,s)e_{k}dr\|_{\dot{H}^{\nu}}^{2}ds\\
&\leq C(H)[(t_{2}-t_{1})^{(2-\nu)\alpha+4H-3}+(t_{2}-t_{1})^{(2-\nu)\alpha+2H-1}]. \tag{3.7}
\end{align*}

In a similar manner, for $\frac{1}{2}<H<1$, we have
\begin{align*}
&\mathbb{E}\|Z(t_{2})-Z(t_{1})\|_{\dot{H}^{\nu}}^{2}\\
&\leq 2\sum\limits_{k=1}^{\infty}\mathbb{E}\|\int_{0}^{t_{1}}\lambda^{1/2}_{k}(K_{t}^{*}(\mathcal{S}_{\alpha}(t_{2}-s)-\mathcal{S}_{\alpha}(t_{1}-s))e_{k})(s)d\beta_{k}(s)\|_{\dot{H}^{\nu}}^{2}\\
&\hspace{2mm}+2\sum\limits_{k=1}^{\infty}\mathbb{E}\|\int_{t_{1}}^{t_{2}}\lambda^{1/2}_{k}(K_{t}^{*}\mathcal{S}_{\alpha}(t_{2}-s)e_{k})(s)d\beta_{k}(s)\|_{\dot{H}^{\nu}}^{2}\\
&= 2\sum\limits_{k=1}^{\infty}\int_{0}^{t_{1}}\mathbb{E}\|\int_{s}^{t}\lambda^{1/2}_{k}(\mathcal{S}_{\alpha}(t_{2}-r)-\mathcal{S}_{\alpha}(t_{1}-r))\frac{\partial K_{H}}{\partial r}(r,s)e_{k}dr\|_{\dot{H}^{\nu}}^{2}ds\\
&\hspace{2mm}+2\sum\limits_{k=1}^{\infty}\int_{t_{1}}^{t_{2}}\|\int_{s}^{t}\lambda^{1/2}_{k}\mathcal{S}_{\alpha}(t_{2}-r)\frac{\partial K_{H}}{\partial r}(r,s)e_{k}dr\|_{\dot{H}^{\nu}}^{2}ds\\
&\leq 2 (t_{2}-t_{1})^{2-(2-\nu)\alpha}\int_{0}^{t_{1}}(\int_{s}^{t}\|\frac{\partial K_{H}}{\partial r}(r,s)\|^{2}dr)(\int_{s}^{t}\mathbb{E}\|\lambda^{1/2}_{k}e_{k}\|^{2}dr)ds\\
&\hspace{2mm}+2\int_{t_{1}}^{t_{2}}(\int_{s}^{t}\|\mathcal{S}_{\alpha}(t_{2}-r)\frac{\partial K_{H}}{\partial r}(r,s)\|_{\dot{H}^{\nu}}^{2}dr)(\int_{s}^{t}\mathbb{E}\|\lambda^{1/2}_{k}e_{k}\|^{2}dr)ds\\
&\leq C(H,Q)[t^{2H}(t_{2}-t_{1})^{2-(2-\nu)\alpha}+(t_{2}-t_{1})^{(2-\nu)\alpha+2H-1}]. \tag{3.8}
\end{align*}

When $H=\frac{1}{2}$, we can deduce that
\begin{align*}
\mathbb{E}\|Z(t_{2})-Z(t_{1})\|_{\dot{H}^{\nu}}^{2}\leq C(H,Q)(t_{2}-t_{1})^{2-(2-\nu)\alpha}. \tag{3.9}
\end{align*}

Therefore, when we set $\gamma=\min\{2-(2-\nu)\alpha,(2-\nu)\alpha+4H-3,(2-\nu)\alpha+2H-1\}>0$ with $\frac{1}{4}<H<1$, taking expectation on both side of (3.5) and combining (3.6)-(3.9) in turn then lead to complete the proof. $\square$

\section{Existence and regularity of mild solution}

In this section, the existence and uniqueness of mild solution to (1.6) will be proved by Banach fixed point theorem. Let $K>0$ be constant to be determined later. We define the following space
\begin{align*}
B_{R}^{T}:=\{u:u\in C([0,T];\dot{H}^{\sigma}),\sup\limits_{t\in[0,T]}\|u(t)\|_{ \dot{H}^{\sigma}}\leq K,~\forall t\in[0,T],\sigma\geq0\},
\end{align*}
where we denote $\dot{H}^{0}:=L^{2}(D)$. The following statement holds.

\textbf{Theorem 4.1.} For $0\leq\nu<2$ and $0<\alpha<1$, then there exists a stopping time $T^{*}>0$ such that (1.6) has a unique mild solution in $L^{2}(\Omega,B_{R}^{T^{*}})$.

\textbf{Proof.} We first define a map $\mathcal{F}:B_{R}^{T}\rightarrow C([0,T];\dot{H}^{\sigma})$ in  the following manner: for any $u\in B_{R}^{T}$,
\begin{align*}
(\mathcal{F}u)(t)&=E_{\alpha}(t)u_{0}+\int_{0}^{t}(t-s)^{\alpha-1}E_{\alpha,\alpha}(t-s)[B(u(s))+f(u(s))]ds\\
&~~~+\int_{0}^{t}(t-s)^{\alpha-1}E_{\alpha,\alpha}(t-s)dB^{H}(s),~t\in [0,T]. \tag{4.1}
\end{align*}

To begin with, we need to show that $\mathcal{F}u\in B_{R}^{T}$ for $u\in B_{R}^{T}$. Making use of Lemma 3.1 and Theorem 3.1 and H\"{o}lder inequality, and based on $\|B(u)\|\leq C\|u\|\|A^{\frac{1}{2}}u\|$, we get
\begin{align*}
\mathbb{E}\|\mathcal{F}u\|_{\dot{H}^{\nu}}^{2}&\leq 3\mathbb{E}\|E_{\alpha}(t)u_{0}\|_{\dot{H}^{\nu}}^{2}+3\mathbb{E}\|\int_{0}^{t}\mathcal{S}_{\alpha}(t-s)[B(u)+f(u)]ds\|_{\dot{H}^{\nu}}^{2}+3\mathbb{E}\|Z(t)\|_{\dot{H}^{\nu}}^{2}\\
&\leq C\mathbb{E}\|u_{0}\|_{\dot{H}^{\nu}}^{2}+Ct^{(2-\nu)\alpha-1}\mathbb{E}(\int_{0}^{t}(\|B(u)\|^{2}+\|f(u)\|^{2})ds)+C(H,Q)t^{\sigma}\\
&\leq C\mathbb{E}\|u_{0}\|_{\dot{H}^{\nu}}^{2}+Ct^{(2-\nu)\alpha}(1+K^{2}+K^{4})+C(H,Q)t^{\sigma}, \tag{4.2}
\end{align*}
which implies that $\mathcal{F}u\in B_{R}^{T}$ as $T>0$ is sufficiently small and $K$ is sufficiently large. By a similar calculation as showing (4.2), we get the continuity of $\mathcal{F}u$.

Given any $u,v\in B_{R}^{T}$, it follows from Lemma 3.1 that
\begin{align*}
\mathbb{E}\|\mathcal{F}u-\mathcal{F}v\|_{\dot{H}^{\nu}}^{2}&\leq 2\mathbb{E}\|\int_{0}^{t}\mathcal{S}_{\alpha}(t-s)[B(u)-B(v)]ds\|_{\dot{H}^{\nu}}^{2}+2\mathbb{E}\|\int_{0}^{t}\mathcal{S}_{\alpha}(t-s)[f(u)-f(v)]ds\|_{\dot{H}^{\nu}}^{2}\\
&\leq Ct^{2\alpha-1}\mathbb{E}(\int_{0}^{t}K^{2}\|u-v\|_{\dot{H}^{\nu}}^{2}ds)+Ct^{2\alpha-1}\mathbb{E}(\int_{0}^{t}\|u-v\|_{\dot{H}^{\nu}}^{2}ds),\tag{4.3}
\end{align*}
which further implies
\begin{align*}
\sup\limits_{t\in[0,T]}\mathbb{E}\|\mathcal{F}u-\mathcal{F}v\|_{\dot{H}^{\nu}}^{2}\leq C(T^{*})^{2\alpha}(1+K^{2})\sup\limits_{t\in[0,T]}\mathbb{E}\|u-v\|_{\dot{H}^{\nu}}^{2}.\tag{4.4}
\end{align*}

Next, let us take $T^{*}$ such that
\begin{align*}
0<C(T^{*})^{2\alpha}(1+K^{2})<1,
\end{align*}
and so that $\mathcal{F}$ is a strict contraction mapping on $B_{R}^{T}$. By the Banach fixed point theorem, there exist a unique fixed point $u\in L^{2}(\Omega,B_{R}^{T^{*}})$, which is a mild solution of (1.6). This completes the proof. $\square$

Our final main result is devoted to the H\"{o}lder regularity of the mild solution and is stated as follows.

\textbf{Theorem 4.2.} For $0\leq\nu<2$, $\frac{1}{4}<H<1$ and $0<\alpha<1$, there exists a unique mild solution $u(t)$ satisfying
\begin{align*}
\mathbb{E}\|u(t_{2})-u(t_{1})\|_{\dot{H}^{\nu}}^{2}< (t_{2}-t_{1})^{\beta}, ~0\leq t_{1}<t_{2}\leq T,
\end{align*}
where $\beta=\min\{\alpha\nu,(2-\nu)\alpha,2-(2-\nu)\alpha,(2-\nu)\alpha+4H-3,(2-\nu)\alpha+2H-1\}>0$.

\textbf{Proof.} From (2.7) we have
\begin{align*}
u(t_{2})-u(t_{1})&=E_{\alpha}(t_{2})u_{0}-E_{\alpha}(t_{1})u_{0}+\int_{0}^{t_{2}}\mathcal{S}_{\alpha}(t_{2}-s)B(u(s))ds-\int_{0}^{t_{1}}\mathcal{S}_{\alpha}(t_{1}-s)B(u(s))ds\\
&\hspace{2mm}+\int_{0}^{t_{2}}\mathcal{S}_{\alpha}(t_{2}-s)f(u(s))ds-\int_{0}^{t_{1}}\mathcal{S}_{\alpha}(t_{1}-s)f(u(s))ds+Z(t_{2})-Z(t_{1})\\
&=:J_{1}+J_{2}+J_{3}+J_{4}, \tag{4.5}
\end{align*}
where we define
\begin{align*}
J_{1}:=E_{\alpha}(t_{2})u_{0}-E_{\alpha}(t_{1})u_{0},~J_{4}:=Z(t_{2})-Z(t_{1}),
\end{align*}
and
\begin{align*}
J_{2}:&=\int_{0}^{t_{2}}\mathcal{S}_{\alpha}(t_{2}-s)B(u(s))ds-\int_{0}^{t_{1}}\mathcal{S}_{\alpha}(t_{1}-s)B(u(s))ds\\
&=\int_{0}^{t_{1}}[\mathcal{S}_{\alpha}(t_{2}-s)-\mathcal{S}_{\alpha}(t_{1}-s)]B(u(s))ds+\int_{t_{1}}^{t_{2}}\mathcal{S}_{\alpha}(t_{2}-s)B(u(s))ds\\
&=:J_{21}+J_{22},
\end{align*}
and
\begin{align*}
J_{3}&:=\int_{0}^{t_{2}}\mathcal{S}_{\alpha}(t_{2}-s)f(u(s))ds-\int_{0}^{t_{1}}\mathcal{S}_{\alpha}(t_{1}-s)f(u(s))ds\\
&=\int_{0}^{t_{1}}[\mathcal{S}_{\alpha}(t_{2}-s)-\mathcal{S}_{\alpha}(t_{1}-s)]f(u(s))ds+\int_{t_{1}}^{t_{2}}\mathcal{S}_{\alpha}(t_{2}-s)f(u(s))ds\\
&=:J_{31}+J_{32}.
\end{align*}

The application of Lemma 2.2 follows that
\begin{align*}
\mathbb{E}\|J_{1}\|_{\dot{H}^{\nu}}^{2}=\mathbb{E}\|E_{\alpha}(t_{2})u_{0}-E_{\alpha}(t_{1})u_{0}\|_{\dot{H}^{\nu}}^{2}\leq (t_{2}-t_{1})^{\alpha\nu}\mathbb{E}\|u_{0}\|^{2}. \tag{4.6}
\end{align*}

Applying the Lemma 3.1 and H\"{o}lder inequality, we get
\begin{align*}
\mathbb{E}\|J_{2}\|_{\dot{H}^{\nu}}^{2}&\leq 2\mathbb{E}\|J_{21}\|_{\dot{H}^{\nu}}^{2}+\mathbb{E}\|J_{22}\|_{\dot{H}^{\nu}}^{2}\\
&\leq2\mathbb{E}(\int_{0}^{t_{1}}\|\mathcal{S}_{\alpha}(t_{2}-s)-\mathcal{S}_{\alpha}(t_{1}-s)\|_{\dot{H}^{\nu}}^{2}ds)(\int_{0}^{t_{1}}\|B(u(s))\|^{2}ds)\\
&\hspace{2mm}+2\mathbb{E}(\int_{t_{1}}^{t_{2}}\|\mathcal{S}_{\alpha}(t_{2}-s)\|_{\dot{H}^{\nu}}^{2}ds)(\int_{t_{1}}^{t_{2}}\|B(u(s))\|^{2}ds)\\
&\leq CK^{4}T^{2}(t_{2}-t_{1})^{2-(2-\nu)\alpha}+CK^{4}(t_{2}-t_{1})^{(2-\nu)\alpha}. \tag{4.7}
\end{align*}
and
\begin{align*}
\mathbb{E}\|J_{3}\|_{\dot{H}^{\nu}}^{2}&\leq 2\mathbb{E}\|J_{31}\|_{\dot{H}^{\nu}}^{2}+\mathbb{E}\|J_{32}\|_{\dot{H}^{\nu}}^{2}\\
&\leq2\mathbb{E}(\int_{0}^{t_{1}}\|\mathcal{S}_{\alpha}(t_{2}-s)-\mathcal{S}_{\alpha}(t_{1}-s)\|_{\dot{H}^{\nu}}^{2}ds)(\int_{0}^{t_{1}}\|f(u(s))\|^{2}ds)\\
&\hspace{2mm}+2\mathbb{E}(\int_{t_{1}}^{t_{2}}\|\mathcal{S}_{\alpha}(t_{2}-s)\|_{\dot{H}^{\nu}}^{2}ds)(\int_{t_{1}}^{t_{2}}\|f(u(s))\|^{2}ds)\\
&\leq C(1+K^{2})T^{2}(t_{2}-t_{1})^{2-(2-\nu)\alpha}+C(1+K^{2})(t_{2}-t_{1})^{(2-\nu)\alpha}.\tag{4.8}
\end{align*}

By Theorem 3.2, we have
\begin{align*}
\mathbb{E}\|J_{4}\|_{\dot{H}^{\nu}}^{2}=\mathbb{E}\|Z(t_{2})-Z(t_{1})\|_{\dot{H}^{\nu}}^{2}\leq C(H,Q)(t_{2}-t_{1})^{\gamma}.\tag{4.9}
\end{align*}

Taking expectation on both side of (4.5) and combining (4.6)-(4.9), the proof of Theorem 4.2 is then completed. $\square$

\section*{Acknowledgements}

Guang-an Zou is supported by National Nature Science Foundation of China (Grant No. 11626085). Guangying Lv is supported by National Nature Science Foundation of China (Grant No. 11771123).


\begin{thebibliography}{00}

\bibitem{Bensoussan-1973} A. Bensoussan, R. Temam, Equations stochastiques du type Navier-Stokes. J. Funct. Anal. 13(2) (1973) 195-222.

\bibitem{Biagini-2008} F. Biagini, Y. Hu, B. {\O}ksendal, T. Zhang, Stochastic calculus for fractional Brownian motion and applications, Springer (2008).

\bibitem{Brzezniak-2013} Z. Brzezniak, T. Caraballo, J.A. Langa, Y. Li, G. Lukaszewiczd, J. Real, Random attractors for stochastic 2D Navier-Stokes
equations in some unbounded domains, J. Differential Equations 255 (2013) 3897-3919.

\bibitem{Da-2002} G. Da Prato, A. Debussche, Two-dimensional Navier-Stokes equations driven by a space-time white noise, J. Funct. Anal. 196(1) (2002) 180-210.

\bibitem{De-2015} P.M. De Carvalho-Neto, P. Gabriela, Mild solutions to the time fractional Navier-Stokes equations in $R^{N}$, J. Differential Equations 259 (2015) 2948-2980.

\bibitem{Duncan-2009} T.E. Duncan, B. Maslowski, B. Pasik-Duncan, Semilinear stochastic equations in a Hilbert space with a fractional Brownian motion, SIAM J. Math. Anal. 40(6) (2009) 2286-2315.

\bibitem{Flandoli-1995}  F. Flandoli, B. Schmalfu{\ss}, Random attractors for the 3D stochastic Navier-Stokes equation with multiplicative
noise, Stoch. Stoch. Rep. 59 (1996) 21-45.

\bibitem{Flandoli-Maslowski-1995} F. Flandoli, B. Maslowski, Ergodicity of the 2-D Navier-Stokes equation under random perturbations, Commun. Math. Phys. 172(1) (1995) 119-141.

\bibitem{Ganji-2010}Z.Z. Ganji, D.D. Ganji, A.D. Ganji, M. Rostamian, Analytical solution of time-fractional Navier-Stokes equation in polar coordinate by homotopy perturbation method, Numer. Meth. Part. D. E. 26(1) (2010) 117-124.

\bibitem{Hairer-2006} M. Hairer, J.C. Mattingly, Ergodicity of the 2D Navier-Stokes equations with degenerate stochastic forcing, Ann. Math. (2006) 993-1032.

\bibitem{Jiang-2012} Y. Jiang, T. Wei, X. Zhou, Stochastic generalized Burgers equations driven by fractional noises, J. Differential Equations 252(2) (2012) 1934-1961.

\bibitem{Kumar-2015} D. Kumar, J. Singh, S. Kumar, A fractional model of Navier-Stokes equation arising in unsteady flow of a viscous fluid, J. Ass. Arab. Univ. Basic Appl. Sci. 17 (2015) 14-19.

\bibitem{Li-2016} X. Li, X. Yang, Y. Zhang, Error estimates of mixed finite element methods for time-fractional Navier-Stokes equations, J. Sci. Comput. 70(2) (2017) 500-515.

\bibitem{Mikulevicius-2005} R. Mikulevicius, B.L. Rozovskii, Global $L_{2}$-solutions of stochastic Navier-Stokes equations, Ann. Probab. 33(1) (2005) 137-176.

\bibitem{Momani-2006} S. Momani, Z. Odibat, Analytical solution of a time-fractional Navier-Stokes equation by Adomian decomposition method, Appl. Math. Comput. 177(2) (2006) 488-494.

\bibitem{Mishura-2008} I.U.S. Mishura, Y. Mishura, Stochastic calculus for fractional Brownian motion and related processes, Springer (2008).

\bibitem{Podlubny-1999} I. Podlubny, Fractional Differential Equations, Academic Press, San Diego, 1999.

\bibitem{Sritharan-2006} S. Sritharan, P. Sundar, Large deviations for the two-dimensional Navier-Stokes equations with multiplicative noise, Stoch. Proc. Appl. 116 (2006) 1636¨C1659.

\bibitem{Srivastava-2006} H.M. Srivastava, J.J. Trujillo, Theory and applications of fractional differential equations, Elsevier (2006).

\bibitem{Taniguchi-2011} T. Taniguchi, The existence of energy solutions to 2-dimensional non-Lipschitz stochastic Navier-Stokes equations in unbounded domains, J. Differential Equations, 251(12) (2011) 3329-3362.

\bibitem{Wang-2015} R. Wang, J. Zhai, T. Zhang, A moderate deviation principle for 2-D stochastic Navier-Stokes equations, J. Differential Equations, 258(10) (2015) 3363-3390.

\bibitem{Wang-Zeng-2010} G. Wang, M. Zeng, B. Guo, Stochastic Burgers' equation driven by fractional Brownian motion, J. Math. Anal. Appl. 371(1) (2010) 210-222.

\bibitem{Xu-2009} T. Xu, T. Zhang, Large deviation principles for 2-D stochastic Navier-Stokes equations driven by L\'{e}vy processes, J. Funct. Anal. 257(5) (2009) 1519-1545.

\bibitem{Zhou-2017} Y. Zhou, L. Peng, On the time-fractional Navier-Stokes equations, Comput. Math. Appl. 73(6) (2017) 874-891.

\bibitem{Zhou-Peng-2017} Y. Zhou, L. Peng, Weak solutions of the time-fractional Navier-Stokes equations and optimal control, Comput. Math. Appl. 73(6) (2017) 1016-1027.

\bibitem{Zhou-Wang-2016} Y. Zhou, J.R. Wang, L. Zhang, Basic theory of fractional differential equations, World Scientific (2016).

\bibitem{Zou-2017} G. Zou,
 B. Wang, Stochastic Burgers equation with fractional derivative driven by multiplicative noise, Comput. Math. Appl. (2017) http://dx.doi.org/10.1016/
    j.camwa.2017.08.023.

\bibitem{Zou-Zhou-2017} G. Zou, Y. Zhou, B. Ahmad, A. Alsaedi, Finite difference/element method for the time-fractional Navier-Stokes equations, (2017) (in revision).







\end{thebibliography}
\end{document}